\documentclass[a4paper,10pt]{article}
\usepackage{amsmath,amsfonts,amssymb,amsthm,enumerate,fullpage}
\usepackage{hyperref,booktabs}
\usepackage[utf8]{inputenc}

\providecommand{\FF}{\mathbb{F}}

\providecommand{\cR}{\mathcal{R}}

\newcommand{\vJ}{\mathbf{J}}
\newcommand{\vj}{\mathbf{j}}

\newcommand{\gauss}[2]{\genfrac{[}{]}{0pt}{}{#1}{#2}}

\newtheorem{theorem}{Theorem}[section]
\newtheorem{lemma}[theorem]{Lemma}
\newtheorem{corollary}[theorem]{Corollary}

\theoremstyle{definition}
\newtheorem{definition}[theorem]{Definition}

\title{The Chromatic Number of the $q$-Kneser Graph for $q \geq 5$}

\author{
Ferdinand Ihringer\footnote{\hbadness=1500 
Department of Mathematics:  Analysis, Logic and Discrete Mathematics, Ghent University, Belgium.
The author is supported by a postdoctoral fellowship of the Research Foundation --- Flanders (FWO).}
}

\begin{document}

\maketitle

\begin{abstract}
  We obtain a new weak Hilton-Milner type result for intersecting families of $k$-spaces in $\FF_q^{2k}$, which improves several known results.
  In particular the chromatic number of the $q$-Kneser graph $qK_{n:k}$ was previously known for $n > 2k$ (except for $n=2k+1$ and $q=2$) or
  $k < q \log q - q$. Our result determines the chromatic number of $qK_{2k:k}$ for $q \geq 5$, so that
  the only remaining open cases are $(n, k) = (2k, k)$ with $q \in \{ 2, 3, 4 \}$ and $(n, k) = (2k+1, k)$ with $q = 2$.
  
\end{abstract}

\section{Introduction}

The Kneser graph $K_{n:k}$ has all $k$-sets of $\{ 1, 2, \ldots, n \}$, $n \geq 2k$, as vertices and two sets are adjacent 
if they are disjoint. The following conjecture due to Kneser \cite{Kneser1955} was shown by Lov\'asz \cite{Lovasz1978}:
 
\begin{theorem}[Lov\'asz (1978)]
  The chromatic number of $K_{n:k}$ is $n-2k+2$.
\end{theorem}

We want to point out that the case $n=2k$ is trivial as $K_{2k:k}$ is bipartite.
A natural generalization is the $q$-analog of the Kneser graph: the $q$-Kneser graph $qK_{n:k}$.
Here we take the $k$-spaces of $\FF_q^{n}$ as vertices and two vertices are adjacent 
if they intersect trivially. Let $\gauss{n}{k}$ denote the number of 
$k$-spaces in $\FF_q^n$. Note that for $0 \leq k \leq n$ we have
\begin{align*}
 \gauss{n}{k} = \prod_{i=1}^k \frac{q^{n-i+1}-1}{q^i-1}.
\end{align*}
Due to two previous results by Blokhuis et al. \cite{Blokhuis2010} for $n > 2k$
and Blokhuis et al. \cite{Blokhuis2012} for $n=2k$, we know the following:

\begin{theorem}[Blokhuis, Brouwer, Chowdhury, Frankl, Mussche, Patk\'{o}s, Sz\H{o}nyi (2010), Blokhuis, Brouwer, Sz\H{o}nyi (2012)]
  If $k \geq 3$ and either $q \geq 3$ and $n \geq 2k+1$, or $q=2$ and $n \geq 2k+2$,
  then the chromatic number of $qK_{n:k}$ is $\gauss{n-k+1}{1}$.
  If either $k < q \log q - q$ or $k \leq 3$, then the chromatic number of $qK_{2k:k}$ 
  is $q^k + q^{k-1}$.
\end{theorem}

We complete this result for $n=2k$ and $q \geq 5$.

\begin{theorem}\label{thm:main_chrom}
  Let $q \geq 5$. Then the chromatic number of $qK_{2k:k}$ is $q^k + q^{k-1}$ 
  for $n = 2k$.
\end{theorem}

The key ingredient of the $n=2k$ proof by Blokhuis et al. \cite{Blokhuis2012} is a weak Hilton-Milner type
result (see \cite{Hilton1967} for the Hilton-Milner theorem for the classical set case, a variation of the 
famous Erd\H{o}s-Ko-Rado theorem \cite{ErdHos1961}). In the following, we use 
projective notation, so we call $1$-spaces \textit{points}, $2$-spaces \textit{lines}, and $(n-1)$-spaces \textit{hyperplanes}.
Call the set of all $k$-spaces on fixed point a \textit{dictator} (also known as point-pencil).
The dual of a dictator consists of all $k$-spaces in a fixed hyperplane.
Due to work by Hsieh \cite{Hsieh1975}, Frankl and Wilson \cite{Frankl1986}, and Godsil and Newman \cite{Newm2004},
we know that the largest independent sets of $qK_{n:k}$ are dictators and, if $n=2k$, duals of dictators, that is
the family of all $k$-spaces in a hyperplane. Blokhuis et al. showed the following:
\begin{theorem}[Blokhuis, Brouwer, Sz\H{o}nyi (2011)]\label{thm:bbs_hm}
  Let $k < q \log q - q$ and let $Y$ be an independent set of $qK_{2k:k}$.
  If $Y$ is not contained in a dictator or its dual, then $|Y| < q^{k(k-1)}/2$.
\end{theorem}
Using a result by Tokushige on cross-intersecting families in vector spaces \cite{Tokushige2013}
and some properties of the spectrum of the Grassmann scheme, we improve this as follows:
\begin{theorem}\label{thm:main_hm}
  Let $q \geq 3$ and $k \geq 5$ and let $Y$ be an independent set of $qK_{2k:k}$.
  If $Y$ is not contained in a dictator or its dual, then
  \begin{align*}
    |Y| < (1 + 3q^{-1}) \gauss{k}{1} \gauss{2k-2}{k-2}.
  \end{align*}
\end{theorem}
As this does not cover $k=4$, we also provide the following:
\begin{theorem}\label{thm:main_hm2}
  Let $q \geq 4$ and let $Y$ be an independent set of $qK_{2k:k}$.
  If $Y$ is not contained in a dictator or its dual, then
  \begin{align*}
    |Y| < 3 \gauss{k}{1} \gauss{2k-2}{k-2}.
  \end{align*}
\end{theorem}
While our second bound is slightly worse than Theorem \ref{thm:bbs_hm} for $k$ and $q$ large, 
it is clearly better for $q$ small compared to $k$. It is easy to construct independent sets
of $k$-spaces of size vaguely $\gauss{k}{1} \gauss{2k-2}{k-2}$, so our result is close
to a proper stability result.

Recently, Cameron-Liebler $k$-space classes (also known as Boolean degree $1$ functions)
received some attention \cite{Blokhuis2018,Filmus2018,Rodgers2018}. In particular, Metsch showed the following \cite{Metsch2017}:

\begin{theorem}[Metsch (2017)]\label{thm:metsch_cl}
  Let $q \geq q_0$ for some universal constant $q_0$ and let $k < q \log q - q - 1$.
  Let $Y$ be a non-trivial Boolean degree $1$ function on $k$-spaces of $\FF_q^{2k}$,
  then $|Y| \geq \frac{q}{5} \gauss{2k-1}{k-1}$.
\end{theorem}
The condition on $k$ comes from Theorem \ref{thm:bbs_hm}, so with Theorem \ref{thm:main_hm} we can improve this to the following:
\begin{theorem}\label{thm:main_cl}
  Let $q \geq q_0$ for some universal constant $q_0$.
  Let $Y$ be a non-trivial Boolean degree $1$ function on $k$-spaces of $\FF_q^{2k}$,
  then $|Y| \geq \frac{q}{5} \gauss{2k-1}{k-1}$.
\end{theorem}
Note that a tedious calculation shows that we can choose $q_0 = 89$ if we follow the argument in \cite{Metsch2017} without optimizing any
of the used constants.
We believe that all Boolean degree $1$ functions for $k > 2$ are trivial, so most likely might be still far from the truth.

Our paper is organized as follows: In Section \ref{sec:grassmann}, we write down some basic background
on the Grassmann scheme, so that we can exploit the spectrum of $qK_{n:k}$. In Section \ref{sec:hm_proof} we prove Theorem \ref{thm:main_hm},
and then we finish our investigation with the mentioned consequences in Section \ref{sec:conseqeunces}
and a conclusion in Section \ref{sec:conclusion}.

\section{The Grassmann Scheme}\label{sec:grassmann}

We summarize some needed notation and results for association schemes in the following.
Delsarte's PhD thesis \cite{Delsarte1973} provides a deeper introduction into the theory 
of combinatorial applications of association schemes.

\begin{definition}
  Let $X$ be a finite set. A $k$-class association scheme is a pair $(X, \cR)$,
 where $\cR = \{ R_0, \ldots R_k \}$ is a set of symmetric binary relations
 on $X$ with the following properties:
 \begin{enumerate}[(a)]
  \item $\{ R_0, \ldots R_k \}$ is a partition of $X \times X$.
  \item $R_0$ is the identity relation.
  \item There are constants $p_{ij}^\ell$ such that for $x, y \in X$ with $(x, y) \in R_\ell$
  there are exactly $p_{ij}^\ell$ elements $z$ with $(x, z) \in R_i$ and $(z, y) \in R_j$.
 \end{enumerate}
\end{definition}
Clearly, $(X, R_i)$ is a $p_0^{ii}$-regular graph. For convenience, we write $v = |X|$. The relation $R_i$ can be described by its adjacency 
matrix $A_i$, so a $(v \times v)$-$0$-$1$-matrix which is the indicator function of $R_i$.
As the $A_i$s are Hermitian and commute, we can diagonalize them simultaneously,
that is their eigenvectors are the same. It is well-known that there are $k+1$ common eigenspaces
$V_0, V_1, \ldots, V_k$ of the $A_i$s. As the $A_i$s are regular, the all-ones vector $\vj$
is an eigenvector and we can assume that $V_0 = \langle \vj \rangle$.
Let $E_i$ denote the orthogonal projection
onto the $i$th eigenspace. We can express the $A_i$s as
\begin{align*}
  A_i = \sum_{j=0}^k P_{ji} E_j.
\end{align*}
Note that $P_{0i} = p_0^{ii}$.

The following stability version of Hoffman's bound for independent sets is surely known
for a long time. Its first application, at least in the context of intersecting families, 
which the author is aware of, is due to Ellis \cite{Ellis2012a}.
We include a proof, limited to the setting of association schemes, to keep this paper mostly self-contained.

\begin{lemma}\label{lem:bnd_rest}
  Let $\chi$ be the characteristic vector of an independent set of $(X, R_i)$.
  Assume that $P_{1i}$ is the smallest eigenvalue of $A_i$ and that $P^-$ is the second smallest
  eigenvalue of $A_i$. Let $E_r$ be the orthogonal projection matrix onto the eigenspaces
  orthogonal to $\langle \vj \rangle$ and the eigenspace of $P_{1i}$. Then
  \begin{align*}
    (P^- - P_{1i}) \chi^T E_r \chi \leq  y \left(-P_{1i} - \frac{P_{0i} - P_{1i}}{v} y \right).
  \end{align*}
\end{lemma}
\begin{proof}
  As $\chi$ is a $0$-$1$-vector, we have
  \begin{align*}
    y = \chi^T \chi = \frac{y^2}{v} + \sum_{i=1}^k \chi^T E_i \chi.%\label{eq:boolean_y}
  \end{align*}
  Hence,
  \begin{align*}
    0 = \chi^T A_i \chi &= \frac{P_{0i}}{v} \chi^T \vJ \chi + \sum_{j=1}^k P_{ji} \chi^T E_j \chi\\
    &\geq \frac{P_{0i}}{v} y^2 + P_{1i} \sum_{j=1}^k \chi^T E_j \chi + (P^- - P_{1i}) \chi^T E_r \chi\\
    &\geq \frac{P_{0i}}{v} y^2 + P_{1i} \left( y - \frac{y^2}{v} \right) + (P^- - P_{1i}) \chi^T E_r \chi.
  \end{align*}
  Rearranging shows the claim.
\end{proof}

The following is surely folklore; see for example the proof of Theorem 2 in \cite{Eisfeld1998} for a mostly identical statement.

\begin{lemma}\label{lem:bnd_rel_i}
  Let $\chi$ be the characteristic vector of a non-empty subset $Y$ of $X$, where $y = |Y|$.
  Let $E_r$ be the orthogonal projection matrix onto the eigenspaces
  orthogonal to $\langle \vj \rangle + V_1$. Let $P^-$ be the smallest
  eigenvalue of $A_i$. Then there exists a $T \in Y$
  such that at least 
  \begin{align*}
    \frac{P_{0i}-P_{1 i}}{v} y + P_{1 i} + (P^- - P_{1 i}) \chi^T E_r \chi / y
  \end{align*}
  elements of $Y$ are in relation $R_i$ to $T$.
\end{lemma}
\begin{proof}
  As in the proof of Lemma \ref{lem:bnd_rest} we obtain
  \begin{align*}
    \chi^T A_i \chi &= \frac{P_{0i}  - P_{1i}}{v} y^2 + P_{1i}y + (P^- - P_{1i}) \chi^T E_r \chi.
  \end{align*}
  Now averaging shows the claim.
\end{proof}

In the \textit{Grassmann scheme} $J_q(n, k)$ the set of all $k$-spaces of $\FF_q^n$ is $X$ and
two subspaces $x$ and $y$ are in relation $R_i$ if their intersection is a subspace of dimension $k-i$.
Clearly, $R_k$ corresponds to adjacency in $qK_{2k:k}$. The eigenvalues $P_{ij}$ of the Grassmann scheme are well-known.
There are two useful formulas, one due to Delsarte \cite{Delsarte1973} and one due to Eisfeld \cite{Eisfeld1999}:
\begin{align}
  P_{ij} &= \sum_{h=0}^i (-1)^{i-h}  q^{hj + \binom{i-h}{2}} 
\gauss{k-j}{h} \gauss{k-h}{i-h} \gauss{n-k-j+h}{h} \label{eq:Pvals1}\\
&= \sum_{h=0}^j (-1)^{j-h}  q^{i(i-j+h) + \binom{j-h}{2}} 
\gauss{j}{h} \gauss{k-h}{i} \gauss{n-k-j+h}{n-k-i}.\label{eq:Pvals}
\end{align}

\section{The Weak Hilton-Milner Theorem}\label{sec:hm_proof}

We rely on the following result by Tokushige. Here a pair $(Y, Z)$, $Y, Z \subseteq X$ is
a \textit{cross-intersecting family} if all elements in $Y$ intersect all elements of $Z$ non-trivially.
Similarly, throughout this section we call an independent set of $qK_{2k:k}$ an \textit{intersecting family}.

\begin{theorem}[Tokushige (2013)]\label{thm:ci_thm}
  Let $(Y, Z)$ be a cross-intersecting family of $qK_{n:k}$. Then
  \begin{align*}
    |Y| \cdot |Z| \leq \gauss{n-1}{k-1}^2.
  \end{align*}
\end{theorem}
For the rest of the section, set $y = (1 + 3q^{-1}) \gauss{k}{1} \gauss{2k-2}{k-2}$.
We also assume that $k > 3$ as the case $k=3$ was taken care of in \cite{Blokhuis2012},
and that $Y$ is not a dictator or the dual of a dictator.

\subsection{Proof of Theorem \ref{thm:main_hm}}\label{sec:gen_sec}

\begin{lemma}\label{lem:quotient}
  Let $\ell$ be a line in $\FF_q^n$. Let $Z$ be a set of $k$-spaces 
  which meet $\ell$ in a fixed point $p$. Set $Z' = \{ \langle z, \ell \rangle/\ell: z \in Z \}$.
  Then $|Z'| \geq |Z| / \gauss{k}{1}$.
\end{lemma}
\begin{proof}
  Let $C$ be a complement of $\ell$ in $\FF_q^n$. For $z \in Z$ we have that
  $\langle z, \ell \rangle$ meets $C$ in a $(k-1)$-space $z'$ as $\dim(C) = n-2$
  and $\dim(\langle z, \ell \rangle) = k+1$.
  There are at most $\gauss{k}{1}$ $k$-spaces $\tilde{z}$ through $p$ in $\langle z', \ell \rangle$.
  Hence, at most $\gauss{k}{1}$ $k$-spaces in $Z$ correspond to the same $(k-1)$-space in $Z'$.
\end{proof}

\begin{lemma}\label{lem:bad_pt_bound}
  Let $k \geq 5$ and $q \geq 3$.
  Let $Y$ be an intersecting family of $qK_{2k:k}$ of size at least $y$, then no point lies in more than $(q^{3-k} \gauss{k}{1}+1) \gauss{2k-2}{k-2}$ elements of $Y$.
\end{lemma}
\begin{proof}
  First we show that no point $p_1$ lies on more than $\gauss{k}{1} \gauss{2k-2}{k-2}$ elements of $Y$.
  This is clear as otherwise there is a $T \in Y$ with $p_1 \notin T$.
  We want to bound the number of $R \in Y$ with $p_1 \in R$.
  We have $\gauss{k}{1}$ choices for one point $p'$ in $R \cap T$
  and then $\gauss{2k-2}{k-2}$ choices for choosing the $k$-space $R$ through $\langle p_1, p' \rangle$.
  Hence, there are at most $\gauss{k}{1} \gauss{2k-2}{k-2}$ elements of $Y$ on $p_1$.

  We continue by showing that if one point $p_1$ lies in at least $(q^{-\alpha} \gauss{k}{1}+1) \gauss{2k-2}{k-2}$ elements of $Y$, then all
  other points lie in at most $(1 + q^{-1}) q^{\alpha+1} \gauss{2k-3}{k-2}$ elements of $Y$ which are not on $p_1$.
  
  Suppose to the contrary that there are two points $p_1$ and $p_2$ such that
  $p_1$ lies in at least $(q^{-1} \gauss{k}{1}+1) \gauss{2k-2}{k-2}$ elements of $Y$.
  Let $Z_1$, respectively, $Z_2$ denote the elements of $Y$ on $p_1$, respectively, $p_2$. Let $\ell$ be $\langle p_1, p_2 \rangle$.
  Let $Z_i' = \{ \langle \ell, z \rangle/\ell : z \in Z_i \text{ and not } \ell \subseteq z \}$ for $i \in \{ 1, 2 \}$.
  Let $x$ be the number of elements of $Y$ containing $\ell$. By Lemma \ref{lem:quotient}, we conclude that
  $|Z_i'| \geq (|Z_i| - x) / \gauss{k}{1}$ and $(Z_1, Z_2)$ is a cross-intersecting family of $(k-1)$-spaces
  in $\FF_q^{2k-2}$. Notice that $x \leq \gauss{2k-2}{k-2}$, so $|Z_1|-x \geq q^{-\alpha} \gauss{k}{1} \gauss{2k-2}{k-2}$. 
  By Theorem \ref{thm:ci_thm}, we obtain
  \begin{align*}
    q^{-\alpha} \gauss{k}{1} \gauss{2k-2}{k-2} \cdot \left( |Z_2| - x \right) \leq \gauss{k}{1}^2 \gauss{2k-3}{k-2}^2.
  \end{align*}
  For $q \geq 3$, this simplifies to
  \begin{align*}
    |Z_2|-x \leq q^\alpha \frac{q^k-1}{q^{2k-2}-1} \gauss{k}{1} \gauss{2k-3}{k-2} \leq (1 + q^{-1}) q^{\alpha+1} \gauss{2k-3}{k-2} =: b.
  \end{align*}
  Let $R \in Y$. As no point on $R$ except for $p_1$ lies in more than $b$ elements of $Y$, $R$ has $\gauss{k}{1}$ points and all elements of $Y$
  meet $R$ in at least one point, we have
  \begin{align*}
    (|Y| - |Z_1|)/\left( \gauss{k}{1} \gauss{2k-2}{k-2} \right) &\leq \frac{q^k-1}{q^{2k-2}-1} \cdot (1 + q^{-1}) q^{\alpha+1} \leq (1 + q^{-1}) q^{3+\alpha-k}.
  \end{align*}
  Suppose that no point $p_1$ lies in at least $q^{-\alpha+1} \gauss{k}{1} \gauss{2k-2}{k-2} > (q^{-\alpha} \gauss{k}{1} + 1) \gauss{2k-2}{k-2}$ elements of $Y$.
  Then
  \begin{align*}
    |Y|/\left( \gauss{k}{1} \gauss{2k-2}{k-2} \right) \leq q^{-\alpha+1} + (1 + q^{-1}) q^{3+\alpha-k}.
  \end{align*}
  Recall that we can assume that $\alpha \geq 1$ as $p_1$ lies in at most $\gauss{k}{1} \gauss{2k-2}{k-2}$ elements of $Y$.
  For $k \geq 5$ and $\alpha \geq 1$, this is less than $1+3q^{-1}$, a contradiction as long as $3+\alpha-k \leq 0$.
  Hence, we can choose $\alpha = k-3$ which shows the assertion.
\end{proof}

\begin{lemma}\label{lem:ci_last_case}
  Let either $k \geq 5$ and $q \geq 3$, or $k \geq 6$ and $q \geq 5$.
  Let $Y$ be an intersecting family of $qK_{2k:k}$ of size at least $y$.
  Then there are no points $p_1$ and $p_2$ such that
  the number of elements of $Y$ on $p_1$ and
  the number of elements of $Y$ on $p_2$ is at least
  \begin{align*}
    (1 - 2q^{3-k} - 2q^{1-k}) y / \gauss{k-1}{1}.
  \end{align*}
\end{lemma}
\begin{proof}
  Let $Z_1$, respectively, $Z_2$ denote the elements of $Y$ on $p_1$, respectively, $p_2$.
  Set $z_1 = |Z_1|$ and $z_2 = |Z_2|$ (we assume $z_1 \geq z_2$).
  Suppose that $z_1, z_2 \geq (1 - 2q^{3-k} - 2q^{1-k}) y / \gauss{k-1}{1}$.
  Set $\ell = \langle p_1, p_2 \rangle$. 
  Let $Z_i' = \{ \langle \ell, z \rangle/\ell : z \in Z_i \text{ and not } \ell \subseteq z \}$ for $i \in \{ 1, 2 \}$.
  As $\ell$ contains at most $\gauss{2k-2}{k-2}$ elements of $Y$, we conclude, using Lemma \ref{lem:quotient}, that
  $|Z_i'| \geq (z_i - \gauss{2k-2}{k-2}) / \gauss{k}{1}$. As $Y$ is an intersecting family, $(Z_1, Z_2)$ is a cross-intersecting family of $(k-1)$-sets
  in $\FF_q^{2k-2}$. By the bound in Theorem \ref{thm:ci_thm}, taking the square root and rearranging, we obtain
  \begin{align*}
     (1 - 2q^{3-k} - 2q^{1-k}) y / \gauss{k-1}{1} - \gauss{2k-2}{k-2} \leq z_2 - \gauss{2k-2}{k-2} \leq \gauss{k}{1} \gauss{2k-3}{k-2}.
  \end{align*}
  By using $y = (1 + 3q^{-1}) \gauss{k}{1} \gauss{2k-2}{k-2}$ and rearranging, we obtain 
  \begin{align*}
    (1 - 2q^{3-k} - 2q^{1-k})(1+3q^{-1}) \leq \frac{q^{k-1}-1}{q-1} \cdot \frac{q^k-1}{q^{2k-2}-1} + \frac{q^{k-1}-1}{q^k-1}.
  \end{align*}
  This is easily verified to be a contradiction under the conditions on $k$ and $q$.
\end{proof}

\begin{lemma} \label{lem:s_cond}
  Let either $k \geq 5$ and $q \geq 3$, or $k \geq 6$ and $q \geq 5$.
  Let $Y$ be an intersecting family of $qK_{2k:k}$ of size at least $y$.
  Let $s$ be the dimension of a smallest subspace meeting all elements of $Y$.
  Then $s \in \{ 1, k \}$.
\end{lemma}
\begin{proof}
  Let $S$ be a subspace meeting all elements of $Y$.
  We suppose that $1 < \dim(S) < k$ and will arrive at a contradiction,
  so suppose that $\dim(S) = k-1$ from now on.
  Let $p_1$ and $p_2$ the points in $S$ which lie on the most elements of $Y$.
  Let $Z_1$, respectively, $Z_2$ denote the elements of $Y$ on $p_1$, respectively, $p_2$.
  Set $z_1 = |Z_1|$ and $z_2 = |Z_2|$ (we assume $z_1 \geq z_2$).
  By Lemma \ref{lem:bad_pt_bound}, $z_1 \leq (q^{3-k} + q^{1-k}) \gauss{k}{1} \gauss{2k-2}{k-2}$.
  Clearly,
  \begin{align*}
    z_2 \geq (y - z_1)/\gauss{k-1}{1} \geq \left( y - (q^{3-k} + q^{1-k}) \gauss{k}{1} \gauss{2k-2}{k-2} \right).
  \end{align*}
  Hence, $z_2 \geq (1 - q^{3-k} - q^{1-k}) y / \gauss{k-1}{1}$. 
  By Lemma \ref{lem:ci_last_case}, this is a contradiction. Hence, $s \in \{ 1, k \}$.
\end{proof}

By duality, we obtain the following:

\begin{corollary}\label{cor:sp_cond}
  Let either $k \geq 5$ and $q \geq 3$, or $k \geq 6$ and $q \geq 5$.
  Let $Y$ be an intersecting family of $qK_{2k:k}$ of size at least $y$.
  Let $s'$ be the dimension of a largest subspace $S$ such that the hyperplanes through $S$ contain all elements of $Y$.
  Then $s' \in \{ 2k-1, k \}$.
\end{corollary}

Hence, in the notation of Lemma \ref{lem:s_cond} and Corollary \ref{cor:sp_cond}, $(s, s') \in \{ (1, 2k-1), (1, k), (k, 2k-1), (k, k) \}$.
If $s = 1$, then all elements of $Y$ lies on a fixed point, so $Y$ is a subset of a dictator.
Similarly, if $s=2k-1$, then all elements of $Y$ lie in a fixed, so is a subset of the dual of a dictator.
Hence, we only need to rule out the case $(s, s') = (k, k)$. 

\begin{lemma}\label{lem:number_in_km1_space}
  Let $Y$ be an intersecting family of $qK_{2k:k}$ of size at least $y$. Then an element in $Y$ meets more than $\gauss{k}{1} \gauss{k-1}{1}$
  elements of $Y$ in a $(k-1)$-space.
\end{lemma}
\begin{proof}
  We assume without loss of generality that $y = |Y|$.
  By Equation \eqref{eq:Pvals1}, we have $P_{0k} = q^{k^2}$, $P_{1k} = -q^{k(k-1)}$ as the smallest eigenvalue and $P_{3k} = -q^{k(k-3) + 3}$
  as the second smallest eigenvalue of $A_k$.
  By Lemma \ref{lem:bnd_rest},
  \begin{align*}
    (P_{3k} - P_{1k}) \chi^T E_r \chi \leq  y(-P_{1k} - \frac{P_{0k} - P_{1k}}{v} y).
  \end{align*}
  We have
  \begin{align*}
    \frac{-P_{1k} - \frac{P_{0k} - P_{1k}}{v} y}{P_{3k} - P_{1k}} &= \frac{q^{k(k-1)} - q^{k^2} (1 + q^{-k}) y / \gauss{2k}{k} }{q^{k(k-1)} ( 1 - q^{-2k+3}) }\\
    &= \frac{ 1 - q^k ( 1 + q^{-k}) y/ \gauss{2k}{k}}{ 1 - q^{-2k+3}}\\
    &= \frac{ 1 - q^k (1 + q^{-k}) (1 + 3q^{-1}) \frac{(q^k-1)^2 (q^{k-1}-1)}{(q-1) (q^{2k}-1)(q^{2k-1}-1)}}{1 - q^{-2k+3}}\\
    &\leq \frac{ 1 - q^{-1} (1 + q^{-k}) (1 + 3q^{-1})(1+q^{-1}) }{1 - q^{-2k+3}} \leq 1 - q^{-1} - 4q^{-2}.
  \end{align*}
  Hence, $\chi^T E_r \chi \leq y(1-q^{-1}-4q^{-2})$.
  We want to apply Lemma \ref{lem:bnd_rel_i} for $i=1$, so we want to show that
  \begin{align*} 
    \frac{P_{01}-P_{11}}{v} y + P_{11} (q^{-1} + 4q^{-2}) + P^- (1 - q^{-1} - 4q^{-2})
  \end{align*}
  is larger than $\gauss{k-1}{1} \gauss{k}{1}$.
  By Equation \eqref{eq:Pvals}, $P^- = -\gauss{k}{1}$, $P_{01} = \gauss{k+1}{1} \gauss{k}{1} - \gauss{k}{1} = q \gauss{k}{1}^2$ and
  $P_{11} = q \gauss{k}{1} \gauss{k-1}{1} - \gauss{k}{1} = q^2 \gauss{k-2}{1} \gauss{k}{1} - 1$.
  Hence, we find
  \begin{align*}
      P_{11} ( q^{-1} + 4q^{-2}) \geq (q+4) \gauss{k}{1} \gauss{k-2}{1} - 2q^{-1}.
  \end{align*}
  and
  \begin{align*}
    P^- (1 - q^{-1} - 4q^{-2}) \leq P^- \leq \gauss{k}{1}.
  \end{align*}
  Hence, as $k \geq 4$, $x$ meets at least $(q+3) \gauss{k}{1} \gauss{k-2}{1}$
  elements of $Y$ in a $(k-1)$-space. It is easily verified that $(q+3) \gauss{k-2}{1} > \gauss{k-1}{1}$.
\end{proof}

\begin{proof}[Proof of Theorem \ref{thm:main_hm}]
  As noted before, we only have to rule out that $(s, s') = (k, k)$ occurs,
  so suppose that $(s, s') = (k, k)$.
  By Lemma \ref{lem:number_in_km1_space}, we can find a $k$-space $R' \in Y$
  which meets more than $\gauss{k}{1} \gauss{k-1}{1}$ elements of $Y$ in a $(k-1)$-space.
  By averaging over the $\gauss{k}{1}$ $(k-1)$-spaces of $R'$, we find a $(k-1)$-space $R$
  that lies in more than $\gauss{k-1}{1}$ elements of $Y$.
  As $s = k$, $R$ is disjoint to one element $T \in Y$. Let $H = \langle R, T \rangle$.
  Set $Z = \{ S \in Y: \dim(S \cap H) = k-1 \}$.
  As there are more than $\gauss{k-1}{1}$ elements through $R$, which are all contained
  in $H$, all elements in $Z$ meet $R$ non-trivially.
  By the dual of Lemma \ref{lem:bad_pt_bound}, $H$ contains at most 
  $(q^{3-k} + q^{1-k}) \gauss{k}{1} \gauss{2k-2}{k-2}$ elements of $Y$.
  Hence,
  \begin{align*}
    |Z| \geq (1 - q^{3-k} - q^{1-k}) y.
  \end{align*}
  By averaging, we find a point $p_1$
  on at least
  \begin{align*}
    z_1 = (1 - q^{3-k} - q^{1-k}) y / \gauss{k-1}{1}
  \end{align*}
  elements of $Z$ and a point $p_2$ in at least
  \begin{align*}
    z_2 = (1 - 2q^{3-k} - 2q^{1-k}) y / \gauss{k-1}{1}
  \end{align*}
  elements of $Z$. By Lemma \ref{lem:ci_last_case}, this is a contradiction.
  Hence, $(s, s') = (k, k)$ does not occur and the proof is complete.
\end{proof}

\subsection{Proof of Theorem \ref{thm:main_hm2}} \label{sec:hm_proof_k_eq_4}

Now $y = 3 \gauss{k}{1} \gauss{2k-2}{k-2}$.
For this case the proof is nearly identical to the proof of Theorem \ref{thm:main_hm}.
Instead of Lemma \ref{lem:bad_pt_bound}, we just use the crude bound that no
point lies on more than $\gauss{k}{1} \gauss{2k-2}{k-2}$ elements of $Y$.
The key difference is that we can replace Lemma \ref{lem:ci_last_case}
with the following.

\begin{lemma}\label{lem:ci_last_case_k_eq_4}
  Let $k \geq 4$.
  Let $Y$ be an intersecting family of $qK_{2k:k}$ of size at least $y$.
  Then there are no points $p_1$ and $p_2$ such that
  the number of elements of $Y$ on $p_1$ or
  the number of elements of $Y$ on $p_2$ is more than
  \begin{align*}
    \frac{1}{3} y / \gauss{k-1}{1}.
  \end{align*}
\end{lemma}
\begin{proof}
  Our setup is as in the proof of Lemma \ref{lem:ci_last_case}, just that this time the resulting inequality is
  \begin{align*}
     \frac{1}{3} y / \gauss{k-1}{1} - \gauss{2k-2}{k-2} \leq z_2 - \gauss{2k-2}{k-2} \leq \gauss{k}{1} \gauss{2k-3}{k-2}.
  \end{align*}
  By using $y = 3 \gauss{k}{1} \gauss{2k-2}{k-2}$ and rearranging, we obtain 
  \begin{align*}
    \frac{q^k-1}{q^{k-1}-1} \leq 1 + \frac{q^k-1}{q^{2k-2}-1} \cdot \frac{q^k-1}{q-1}.
  \end{align*}
  This is a contradiction. The assertion follows.
\end{proof}
From here on it is easy to copy the steps which we took for the proof of Theorem \ref{thm:main_hm},
replacing Lemma \ref{lem:bad_pt_bound} and Lemma \ref{lem:ci_last_case} accordingly.

\section{The Chromatic Number}\label{sec:conseqeunces}

In \cite[p. 192]{Blokhuis2012} it was established that if $q^k+q^{k-1}$ is not the chromatic number
and $f$ is the size of the largest independent set which is not contained in a dictator or its dual,
then
\begin{align}
  (q^k - q^{k-1}) \gauss{2k-1}{k-1} q^{k-1} - f \gauss{k}{1} \gauss{k+1}{1} < \epsilon \left( 2f - \gauss{2k-1}{k-1} \right). \label{eq:chrom}
\end{align}
for some $\epsilon > 0$. By Theorem \ref{thm:main_hm}, 
\begin{align*}
  2f / \gauss{2k-1}{k-1} &\leq 2 (1 + 3q^{-1}) \cdot \frac{q^k-1}{q-1} \cdot \frac{q^{k-1}-1}{q^{2k-1}-1}.
\end{align*}
For $q \geq 5$ this is easily verified to be less than $1$ and therefore the right hand side of Equation \eqref{eq:chrom}
is negative. Similarly,
\begin{align*}
  \frac{f \gauss{k}{1} \gauss{k+1}{1}}{(q^k - q^{k-1}) \gauss{2k-1}{k-1} q^{k-1}} &\leq \frac{q^{k-1}-1}{q^{2k-1}-1} \cdot \frac{(q^k-1)^2(q^{k+1}-1)}{(q^{2k-1}-q^{k-1})(q-1)^3} < 1,
\end{align*}
so the left hand side of Equation \eqref{eq:chrom} is positive.
As this is a contradiction, we have shown Theorem \ref{thm:main_chrom}.
Note that \cite[Proposition 5.1]{Blokhuis2012} gives a characterization of the case of equality.

Theorem \ref{thm:main_cl} is a simple consequence of replacing Theorem \ref{thm:bbs_hm} with Theorem \ref{thm:main_hm}
in the proof of Theorem \ref{thm:metsch_cl}. See \cite{Metsch2017} for details.

\section{Future Work}\label{sec:conclusion}

Clearly, the most urgent open cases are the determination of the chromatic number $qK_{2k:k}$
for $q = 2, 3, 4$. 
For $q=3, 4$ it is sufficient to obtain slightly better stability type results which show 
$f \leq \gauss{k}{1} \gauss{2k-2}{k-2}$ as then $2f < \gauss{2k-1}{k-1}$.
For $q=2$ the current approach of determining the chromatic number cannot work as there are examples
very close in size to a dictator and its dual.

Classical polar spaces are the geometries induced by non-degenerate sesquilinear forms onto $\FF_q^n$.
There are currently barely any stability results known for intersecting families of maximals of finite classical polar spaces in literature
and there is an interesting diversity of largest families \cite{Pepe2011}, similar to $qK_{2k:k}$. Results for cross-intersecting families are known
for finite classical polar spaces \cite{Ihringer2015}, so it might be feasible to determine their chromatic number in a
similar fashion.

\section*{Acknowledgements}
I would like to thank the anomymous referee and Jozefien D'haeseleer for their comments
on an earlier draft of this document.

% \bibliographystyle{plain}
% \bibliography{../../jabref}

\begin{thebibliography}{10}

\bibitem{Blokhuis2010}
A.~Blokhuis, A.~E. Brouwer, A.~Chowdhury, P.~Frankl, T.~Mussche, B.~Patk{\'o}s,
  and T.~Sz{\H{o}}nyi.
\newblock A {H}ilton-{M}ilner theorem for vector spaces.
\newblock {\em Electron. J. Combin.}, 17(1):Research Paper 71, 12, 2010.

\bibitem{Blokhuis2012}
A.~Blokhuis, A.~E. Brouwer, and T.~Sz{\H o}nyi.
\newblock On the chromatic number of {$q$}-{K}neser graphs.
\newblock {\em Des. Codes Cryptogr.}, 65(3):187--197, 2012.

\bibitem{Blokhuis2018}
A.~Blokhuis, M.~De~Boeck, and J.~D'haeseleer.
\newblock {C}ameron-{L}iebler sets of $k$-spaces in {PG}$(n,q)$.
\newblock {\em Des. Codes Cryptogr.}, 2018.

\bibitem{Delsarte1973}
P.~Delsarte.
\newblock An algebraic approach to the association schemes of coding theory.
\newblock {\em Philips Res. Rep. Suppl.}, (10):vi+97, 1973.

\bibitem{Eisfeld1998}
J.~Eisfeld.
\newblock Subsets of association schemes corresponding to eigenvectors of the
  {B}ose-{M}esner algebra.
\newblock {\em Bull. Belg. Math. Soc. Simon Stevin}, 5(2-3):265--274, 1998.
\newblock Finite geometry and combinatorics (Deinze, 1997).

\bibitem{Eisfeld1999}
J.~Eisfeld.
\newblock The eigenspaces of the {B}ose-{M}esner algebras of the association
  schemes corresponding to projective spaces and polar spaces.
\newblock {\em Des. Codes Cryptogr.}, 17(1-3):129--150, 1999.

\bibitem{Ellis2012a}
D.~Ellis.
\newblock A proof of the {C}ameron-{K}u conjecture.
\newblock {\em J. Lond. Math. Soc. (2)}, 85(1):165--190, 2012.

\bibitem{ErdHos1961}
P.~Erd{\H{o}}s, C.~Ko, and R.~Rado.
\newblock Intersection theorems for systems of finite sets.
\newblock {\em Quart. J. Math. Oxford Ser. (2)}, 12:313--320, 1961.

\bibitem{Filmus2018}
Y.~Filmus and F.~Ihringer.
\newblock Boolean degree 1 functions on some classical association schemes.
\newblock {\em J. Combin. Theory Ser. A}, 162:241--270, 2019.

\bibitem{Frankl1986}
P.~Frankl and R.~M. Wilson.
\newblock The {E}rd{\H o}s-{K}o-{R}ado theorem for vector spaces.
\newblock {\em J. Combin. Theory Ser. A}, 43(2):228--236, 1986.

\bibitem{Hilton1967}
A.~J.~W. Hilton and E.~C. Milner.
\newblock Some intersection theorems for systems of finite sets.
\newblock {\em Quart. J. Math. Oxford Ser. (2)}, 18:369--384, 1967.

\bibitem{Hsieh1975}
W.~N.~Hsieh.
\newblock Intersection theorems for systems of finite vector spaces.
\newblock {\em Discrete Math.}, 12:1--16, 1975.

\bibitem{Ihringer2015}
F.~Ihringer.
\newblock Cross-intersecting {E}rd{\H o}s-{K}o-{R}ado sets in finite classical polar spaces.
\newblock {\em  Electron. J. Combin.}, 22(2), 2015.

\bibitem{Kneser1955}
M.~Kneser.
\newblock Aufgabe 300.
\newblock In {\em Jber. Deutsch. Math.-Verein. 58}. 1955.

\bibitem{Lovasz1978}
L.~Lov\'asz.
\newblock Kneser's conjecture, chromatic number, and homotopy.
\newblock {\em J. Combin. Theory Ser. A}, 25(3):319--324, 1978.

\bibitem{Metsch2017}
K~.Metsch.
\newblock A gap result for {C}ameron-{L}iebler {$k$}-classes.
\newblock {\em Discrete Math.}, 340(6):1311--1318, 2017.

\bibitem{Newm2004}
M.~W. Newman.
\newblock {\em Independent Sets and Eigenspaces}.
\newblock PhD thesis, University of Waterloo, Waterloo, Canada, 2004.

\bibitem{Rodgers2018}
M.~Rodgers, L.~Storme, and A.~Vansweevelt.
\newblock The {C}ameron-{L}iebler $k$-classes in ${PG}(2k+1, q)$.
\newblock {\em Combinatorica}, 38(3):739--757.

\bibitem{Pepe2011}
V.~Pepe, L.~Storme, and F.~Vanhove.
\newblock Theorems of {E}rd?s-{K}o-{R}ado type in polar spaces.
\newblock {\em J. Combin. Theory Ser. A}, 118(4):1291--1312.

\bibitem{Tokushige2013}
N.~Tokushige.
\newblock The eigenvalue method for cross {$t$}-intersecting families.
\newblock {\em J. Algebraic Combin.}, 38(3):653--662, 2013.

\end{thebibliography}

\end{document}